\documentclass[12pt]{article}
\usepackage{amsfonts}
\begin{document}
\title{\bf Absolute continuity of the representing measures of the Dunkl
 intertwining operator and of its dual and applications
\author{ Khalifa TRIM\`ECHE\\
{\footnotesize{\em Department of Mathematics}}\\
{\footnotesize{\em Faculty of sciences of Tunis- CAMPUS-}}\\
{\footnotesize{\em 1060.Tunis. Tunisia.}} }}
\date{}
\maketitle
\begin{abstract}
In this paper  we prove the absolute continuity of the
representing measures of the Dunkl intertwining operator and of
its  dual. Next we present some applications of this result.
\end{abstract}
{\small{\bf Key word : }} {\small Dunkl intertwining operator and
its dual. Absolute continuity of the representing measures}. 

\bigskip
\noindent{\bf{Introduction}}\\ \hspace*{5mm} \hspace*{5mm} We
consider the differential-difference operators  on $\mathbb{R}^{d}$
introduced by C.F.Dunkl in [4] and called Dunkl operators in the
literature . These operators are very important in pure Mathematics
and in Physics.They  provide a useful tool in the study of special
functions with root systems (see [3] [8]), and they are closely
related to certain representations of degenerate affine Hecke
algebras [2][16], moreover the commutative algebra generated by
these operators has been used in the study of certain
 exactly
solvable
 models of quantum mechanics, namely the
 Calogero-Sutherland-Moser models,
which deal with systems of identical particles in a one
dimensional spaces
 (see [10] [13] [14]).\\

\hspace*{5mm}C.F.Dunkl has proved in [6] that there exists a
unique isomorphism $V_k$ from the space of homogeneous polynomial
${\cal P}_n$ on $\mathbb{R}^d$ of degree $n$ onto itself
satisfying the transmutation relations $$ T_j V_k  =  V_k
\frac{\partial}{\partial x_j}, \quad j = 1, ... , d,\eqno{(1)}$$
$$V_k(1) =  1.\eqno{(2)} $$ This operator is called Dunkl
intertwining operator. Next K.Trim\`eche has extended this
operator to an isomorphism from ${\cal E}(\mathbb{R}^d)$
 (the space of
$C^{\infty}$-functions on $\mathbb{R}^{d}$) onto itself satisfying
the relations (1) and (2) (see [23]).\\

\hspace*{5mm} The operator $V_k$ possesses the integral
representation $$
  V_k (f)(x) = \displaystyle\int_{\mathbb{R}^d}f(y)d\mu_{x}(y), \;  f \in {\cal
  E}(\mathbb{R}^d),\eqno{(3)}$$
where $\mu_x$ is a probability measure on $\mathbb{R}^d$ with
support in the closed ball $B(o,||x||)$ of center $o$ and radius
$||x||$. (See [17][23]).\\

 We have studied in [23] the transposed
operator $^{t}V_k$ of the operator $V_k$. It has the integral
representation $$^{t} V_k (f)(y) =
\displaystyle\int_{\mathbb{R}^d}f(x)d\nu_{y}(x),\eqno{(4)} $$
where $\nu_y$ is a positive measure on $\mathbb{R}^d$ with support
in the set $\{x \in \mathbb{R}^d / ||x|| \geq ||y||\}$ and $f$ in
$D(\mathbb{R}^d)$ (the space of $C^{\infty}$-functions on
$\mathbb{R}^d$ with compact support). This operator is called the
dual Dunkl intertwining operator.

We have proved in [23] that the operator $^{t}V_k$ is an isomorphism
from $D(\mathbb{R}^d)$ onto itself, satisfying the transmutation
relations $$ \forall \, y \in \mathbb{R}^d, \; ^{t}V_k(T_j f)(y) =
\frac{\partial}{\partial y_j}\; ^{t}V_k (f)(y), \;
j=1,...,d,\eqno{(5)}$$ In this paper we prove that the measure
$\mu_x$ given by (3), is absolutely continuous with respect to the
Lebesgue measure on $\mathbb{R}^d$. More precisely for all
continuous function $g$ on $\mathbb{R}^d$, we have
$$\forall\, x \in \mathbb{R}^d, \quad \omega_k(x)V_k(g)(x) =
\displaystyle\int_{\mathbb{R}^d}\mathcal{K}^{o}(x,y)
g(y)dy,\eqno{(6)}$$ and  $$\forall\, x \in
\mathbb{R}^d_{\mbox{reg}}, \quad V_k(g)(x) =
\displaystyle\int_{\mathbb{R}^d}\mathcal{K}(x,y) g(y)dy,\eqno{(7)}
$$ where $\mathcal{K}^{o}(x,.)$ is a positive integrable function
on $\mathbb{R}^d$ with respect to the Lebesgue measure and with
support in $\{y \in \mathbb{R}^d/ \|y\| \leq \|x\|\}$, and
$\mathcal{K}(x,y)$ the function given by

 $$\forall\, x \in
\mathbb{R}^d_{\mbox{reg}}, \forall \, y \in \mathbb{R}^d, \;
\mathcal{K}(x,y) = \omega_k^{-1}(x)\mathcal{K}^{o}(x,y).
\eqno{(8)}$$

 Next we
establish that for  all $y \in \mathbb{R}^d$ the measure $\nu_y$
given by (4), is absolutely continuous with respect to the measure
$\omega_k(x)dx$ on $\mathbb{R}^d$, with $\omega_k$ a positive
weight function on $\mathbb{R}^d$ which will be given in the
following section. More precisely   for all continuous function
$f$ on $\mathbb{R}^d$ with compact support, we have $$ \forall\; y
\in \mathbb{R}^d,\;\;{}^tV_k(f)(y) =
\displaystyle\int_{\mathbb{R}^d}\mathcal{K}(x,y) f(x)
\omega_k(x)dx,\eqno{(9)}$$ where ${\cal K}(.,y)$ is the function
given by the relation (8). It is locally integrable
 on $\mathbb{R}^d$ with support in $\{x \in
\mathbb{R}^d / \|x\| \geq \|y\|\}$.\\

We present some applications of the relations (6),(7), in
particular we prove that  the Dunkl kernel $K(-ix,z)$ satisfies
$$\forall\, x \in \mathbb{R}^d, \; \lim_{\|z\| \rightarrow +
\infty}\{\omega_{k}(x)K(-ix, z)\} = 0,\eqno{(10)}$$ and
 $$\forall \,x \in
\mathbb{R}^d_{reg}, \;\lim_{\|z\| \rightarrow +
 \infty}K(-ix, z) =
0.\eqno{(11)}$$
 Also we give a simple proof of the main result of
[25].\\ \hspace*{5mm} Finally we remark that in personal
communications sent to C.F.Dunkl, M.F.E.de Jeu and M.R\"{o}sler
after the summer of Year 2000, we have conjectured that the measures
$\mu_x$ and $\nu_y$ are absolutely continuous, and we have tried to
solve this conjecture. Next M.F.E.de Jeu and M.R\"{o}sler have also
conjectured in [12] that the measure $\mu_x$ is absolutely
continuous.
\section{The eigenfunction of the Dunkl operators }
In this section we collect some notations and results on Dunkl
operators and the Dunkl kernel (see [5],[6],[9],[11]).
\subsection{Reflection groups, root systems and
multiplicity functions}
We consider ${\mathbb{R}}^{d}$ with the euclidean  scalar product
$\langle . , .\rangle$ and $||x|| = \sqrt{\langle x , x \rangle}$.
On ${\mathbb{C}}^{d},\;||.||$ denotes also the standard Hermitian
norm, while $\langle z, w \rangle = \sum_{j =1 }^{d} z_{j}
\overline{w}_{j} $.\\ \hspace*{5mm}For $\alpha \in
\mathbb{R}^{d}\backslash\{0\} $, let $\sigma_{\alpha}$ be the
reflection in the
 hyperplan $H_{\alpha} \subset \mathbb{R}^{d} $ orthogonal to $\alpha$, i.e.
$$ \sigma_{\alpha}(x) = x - 2 \frac{\langle\alpha, x\rangle}
{||\alpha||^{2}} \alpha.\eqno{(1.1)}$$
\hspace*{5mm} A finite set  $R \subset
\mathbb{R}^{d}\backslash\{0\} $
 is called a root system if
$R \cap {\mathbb{R}}^{d}.\alpha = \{\alpha, - \alpha\}$ and
$\sigma_{\alpha} R = R$ for all $\alpha \in R$. We assume that it
is normalized by $||\alpha||^2 = 2$ for all $\alpha \in R$. For a
given root system $R$
 the reflections $\sigma_{\alpha}, \alpha \in R $,
generate a finite group $W \subset O(d)$,  the reflection group
associated with $R$. All reflections in W correspond to suitable
pairs of roots. For a given $\beta \in \mathbb{R}^{d}_{reg}=
\mathbb{R}^{d} \backslash \cup_{\alpha \in R}H_{\alpha}$, we fix
the positive subsystem $R_{+} = \{\alpha \in R\; / \langle\alpha,
\beta\rangle > 0 \}$, then for each $\alpha \in R$ either $\alpha
\in R_{+}$ or
 $- \alpha \in R_{+}$.\\
  \hspace*{5mm} A function $k :R
\longrightarrow \mathbb{C}$ on a root system R is called a
multiplicity function if it is invariant under the action of the
associated reflection group W. If one regards k as a function on
the corresponding reflections, this means that k is constant on
the conjugacy classes of reflections in W.  For abbreviation, we
introduce the index
$$ \gamma = \gamma(R) = \sum_{\alpha \in R_{+}
}k(\alpha).\eqno{(1.2)}$$
Moreover, let $\omega_{k}$ denotes the weight function
 $$
\omega_{k}(x) = \prod_{\alpha \in R_{+}
}|\langle\alpha,x\rangle|^{2k(\alpha)},\eqno{(1.3)}$$
which is $W-$invariant and homogeneous of degree $2 \gamma$.\\\\
\hspace*{5mm} For $d=1$ and $W= \mathbb{Z}_2$, the multiplicity
function $k$ is a single parameter denoted $\gamma > 0$ and $$
\forall \, x \in\mathbb{ R}, \; \omega_{k}(x) = |x|^{2\gamma}.
\eqno{(1.4)}$$ \hspace*{5mm} We introduce the Mehta-type constant
$$ c_{k} =
(\displaystyle\int_{\mathbb{R}^{d}}\exp(-||x||^{2})\omega_{k}
(x)\;dx)^{-1},\eqno{(1.5)}$$
which is known for all Coxeter groups $W$.(See [4][8][15]).\\
 \hspace*{5mm}
For an integrable function on $\mathbb{R}^{d}$
with respect to the
measure $\omega_{k}(x)\;dx$ we have the relation
$$ \displaystyle\int_{\mathbb{R}^{d}} f(x)\omega_{k}(x)\;dx =
\displaystyle\int_{0}^{+\infty}(\displaystyle\int_{S^{d-1}}f(r
\beta) \omega_{k}( \beta)\;d\sigma(\beta))
r^{2\gamma+ d-1}\;dr,
\eqno{(1.6)}$$
where $d\sigma$ is the normalized surface measure on the unit
sphere
 $S^{d-1}$ of $\mathbb{R}^{d}$.\\
In particular if f is radial (i.e. SO(d)-invariant ), then there
exists a function F on $[0, + \infty[$, such that $f(x) = F(||x||)
= F(r)$, with $||x|| = r$, and the relation (1.6) takes the form
$$ \displaystyle\int_{\mathbb{R}^{d}} f(x)\omega_{k}(x)\;dx = d_k
\displaystyle\int_{0}^{+\infty} F(r) r^{2\gamma + d
-1}\;dr,\eqno{(1.7)}$$
where
$$ d_{k} = \displaystyle\int_{S^{d-1}}\omega_{k}(\beta)\;d\sigma(\beta) =
\frac{2}{c_k \Gamma (\gamma + 2 d)}.\eqno{(1.8)}$$

\subsection{Dunkl operators and Dunkl kernel}
The Dunkl operators $T_{j}\; j\; = 1\;, ...,\; d $, on
$\mathbb{R}^{d}$ associated with the
 finite reflection group W and multiplicity function k are given for a function $f$ of
  class $C^1$ on $\mathbb{R}^d$ by
$$ T_{j} f(x) = \frac{\partial}{\partial x_{j}} f(x) +
\sum_{\alpha \in R_{+}}k(\alpha) \alpha_{j} \frac{f(x) -
f(\sigma_{\alpha}(x))}{\langle\alpha,x\rangle}.\eqno{(1.9)}$$
In the case $k = 0$, the $T_{j}, \, j = 1, ... , d,$ reduce to the
corresponding partial derivatives. In this paper, we will assume
throughout that $k \geq 0$ and $\gamma \geq 0$.\\ \hspace*{5mm} For
$f$ of class $C^1$ on $ \mathbb{R}^{d}$ with compact support and $g$
of class $C^{1}$ on $\mathbb{R}^{d}$ we have
$$ \displaystyle\int_{\mathbb{R}^{d}} T_{j}f(x) g(x)\omega_{k}(x)\;dx = -
\displaystyle\int_{\mathbb{R}^{d}} T_{j}g(x)
 f(x)\omega_{k}(x)\;dx.\eqno{(1.10)}$$
\hspace*{5mm} For $y \in \mathbb{R}^{d} $,  the system $$ \left\{
\begin{array}{crll}
T_{j}u(x,y) &=& y_{j} u(x,y),& j = 1, ..., d,\\\\ u(0,y) &=& 1,

\end{array}
\right. \eqno{(1.11)}$$ admits a unique analytic solution on
$\mathbb{R}^{d}$, which will be denoted $K(x,y)$ and called the
Dunkl kernel.\\ This kernel has a unique holomorphic extension to
${\mathbb{C}}^{d} \times {\mathbb{C}}^{d}$.\vspace*{3mm}\\
\noindent{\bf{Examples.1.1}}\\ \hspace*{5mm} 1) If $d = 1$ and $W
={ \mathbb{Z}}_{2}$, the Dunkl kernel is given by
$$ K(z,t) = j_{\gamma - \frac{1}{2}}(izt) + \frac{zt}{2 \gamma +
1}
 j_{\gamma + \frac{1}{2}}(izt), \quad z, \; t \in \mathbb{C},\eqno{(1.12)}$$
where for $\alpha \geq \frac{-1}{2}$, $j_{\alpha}$ is the
normalized Bessel function
 defined by
$$ j_{\alpha}(z) = 2^{\alpha} \Gamma(\alpha + 1)
\frac{J_{\alpha}(z)}{z^{\alpha}} =  \Gamma(\alpha + 1) \sum_{n =
0}^{\infty}\frac{(-1)^{n}(\frac{z}{2})^{2 n} } {n!  \Gamma(\alpha
+ n + 1)},\eqno{(1.13)}$$
with $J_{\alpha}$ is the Bessel function of first kind and index
$\alpha$. (See [6]).\\ \hspace*{5mm}2) The Dunkl kernel of index
$\gamma = \displaystyle \sum_{l = 1}^{d}\alpha_l, \; \alpha_l >
0$, associated with the reflection group $\mathbb{Z}_2 \times ...
\times \mathbb{Z}_2$ on $\mathbb{R}^d$ is given for all $x, y \in
\mathbb{R}^d$ by
 $$K(x,y) =
\displaystyle \prod_{l = 1}^{d}K(x_l, y_l),\eqno{(1.14)} $$ where
$K(x_l, y_l)$ is the function defined by (1.12).\\

 \hspace*{5mm} The Dunkl kernel possesses the following
 properties.\\
\hspace*{5mm}  i) For $z, t \in \mathbb{C}^{d}$, we have$ K(z,t) =
K(t,z)$; $K(z,0) = 1$ and $K(\lambda z,t) = K(z, \lambda t),$ for
 $ \lambda \in \mathbb{C}.$\\\hspace*{5mm}  ii)
For all $x, y \in \mathbb{R}^{d}$ we have
$$ |K(i x , y)| \leq 1,\eqno{(1.15)}$$
\label{P1.2}
 \hspace*{5mm}$iii)$ The function $K(x,z)$ admits for all
 $x \in \mathbb{R}^{d}$ and $z \in \mathbb{C}^{d}$ the following Laplace type integral representation
$$ K(x,z) = \displaystyle\int_{\mathbb{R}^d} e^{<y,z>}
d\mu_{x}(y),\eqno{(1.16)}$$
where $\mu_{x}$ is the  measure given by the relation (3)
satisfying\\ $\bullet\quad supp\; \mu_{x} \cap \{ y \in
\mathbb{R}^{d} / ||y||=||x||\} \neq \O.$ \hfill(1.17)\\
$\bullet\quad$ For each $r
> 0, w \in W$ and each Borel set $E \subset \mathbb{R}^d$ we have
$$\mu_{rx}(E) = \mu_x(r^{-1}E),\mbox{ and } \mu_{wx}(E) =
\mu_x(w^{-1}E),\eqno{(1.18)}$$ (See [17]).\vspace*{3mm}\\
 \noindent{\bf{Examples 1.2}}\\
\hspace*{5mm}1) When $d=1$ and $W = \mathbb{Z}_2$, for all $x \in
\mathbb{R} \backslash \{0\}$ and $z \in \mathbb{C}$ the relation
(1.16) is of the form $$
   K(x,z) = \frac{\Gamma(\gamma+\frac{1}{2})}{\sqrt{\pi}\Gamma(\gamma)}
|x|^{-2\gamma}\displaystyle\int_{-|x|}^{|x|}(|x|-y)^{\gamma}(|x|+y)^{\gamma-1}e^{yz}dy.
\eqno{(1.19)}$$
Then in this case for all $x \in \mathbb{R} \backslash  \{0\}$ the
measure $\mu_{x}$ is given by $d\mu_{x}(y) = {\cal K}(x,y)dy $
with $$
   {\cal K}(x,y) =  \frac{\Gamma(\gamma+\frac{1}{2})}{\sqrt{\pi}\Gamma(\gamma)}
|x|^{-2\gamma}(|x|-y)^{\gamma}(|x|+y)^{\gamma-1}1_{]-|x|,|x|[}(y),
\eqno{(1.20)} $$ where $1_{]-|x|,|x|[}$ is the characteristic
function of the interval $]-|x|,|x|[$.\\ \hspace*{5mm}2) The Dunkl
kernel of index $\gamma = \displaystyle \sum_{l = 1}^{d}\alpha_l,
\; \alpha_l > 0$, associated with the reflection group
$\mathbb{Z}_2 \times ... \times \mathbb{Z}_2$ on $\mathbb{R}^d$,
possesses
 for all $x \in \mathbb{R}^{d}_{reg} = \mathbb{R}^{d}\backslash \displaystyle\bigcup_{l=1}^{d}H_l $, with
  $H_l = \{x \in \mathbb{R}^d/ \, x_l = 0\}$, and $z \in \mathbb{C}^d$,
   the integral representation
$$K(x,z) = \displaystyle\int_{\mathbb{R}^d} {\cal K}(x,y)e^{\langle y, z
\rangle} dy,\eqno{(1.21)}$$ where
 $${\cal K}(x,y) =
\displaystyle \prod_{l = 1}^{d}{\cal K}(x_l, y_l),\eqno{(1.22)}
,$$ with ${\cal K}(x_l, y_l)$ given by the relation  (1.20).
\section{The Dunkl intertwining operator and its dual}
{\bf{Notations}}.  We denote by   $C(\mathbb{R}^{d}) (resp.
\;C_{c} (\mathbb{R}^{d}))$\, the space of continuous functions on
$\mathbb{R}^{d}$ (resp. with compact support).\vspace*{3mm}\\
 \hspace*{5mm}The Dunkl
intertwining operator $V_k$ is defined on $C(\mathbb{R}^{d})$ by
$$
 \forall\; x \in \mathbb{R}^d ,\quad V_k f(x) = \displaystyle\int_{\mathbb{R}^d}
f(y)d\mu_{x}(y), \eqno{(2.1)}$$
where $\mu_{x}$ is the measure given by the relation (3). (See
[17][23]p.364-366).\\ It possesses many properties in particular we
have
\\ \hspace*{5mm} i) For all $g$  in $C(\mathbb{R}^d)$ the function
$V_k (g)$ belongs to $C(\mathbb{R}^d)$. Moreover for all $x \in
\mathbb{R}^d$ in the closed ball $\overline{B}(o, a)$ of center
$o$ and radius $a
> 0$, we have
$$
  |V_k (g)(x)| \leq \sup_{y \in \overline{B}(o,a)}|g(y)|.\eqno{(2.2)}$$
\hspace*{5mm}ii) We have $$ V_k (g)(o) =
g(o).\eqno{(2.3)}$$\hspace*{4mm}
 iii) We
have $$\forall x \in \mathbb{R}^{d}, \; \; \forall z \in
\mathbb{C}^{d}, \; \; K(x,z) = V_k (e^{\langle .,z\rangle})(x).
\eqno{(2.4)}$$
   \hspace*{5mm}
 The operator
 $^{t}V_k$ satisfying for
$f $ in $C_c(\mathbb{R}^d)$ and $g$ in $C(\mathbb{R}^d)$ the
relation
$$ \displaystyle\int_{\mathbb{R}^d}\, ^{t}V_k(f)(y)g(y)dy = \displaystyle\int_{\mathbb{R}^d}
V_k (g)(x)f(x)\omega_k (x)dx,\eqno{(2.5)} $$
 is given by
$$ \forall\; y \in \mathbb{R}^d,\;\;^{t}V_k(f)(y) =
\displaystyle\int_{\mathbb{R}^d} f(x)d\nu_{y}(x), \eqno{(2.6)} $$
where $\nu_{y}$ is a positive measure on the $\sigma$-algebra of
$\mathbb{R}^d$, satisfying
\\$\bullet$ $supp \, \nu_y \subset \{x \in \mathbb{R}^d / ||x|| \geq
||y||\}.$\\
 $\bullet$ For each $a > 0$, $w \in W$ and each Borel
set $E \subset \mathbb{R}^d$ we have $$\nu_{ay}(E) =
a^{2\gamma}\nu_y(a^{-1}E) \; \mbox{and} \; \nu_{wy}(E) =
\nu_y(w^{-1}E).\eqno{(2.7)}$$ The operator $^{t}V_k$ is called the
dual Dunkl intertwining operator. (See [23]p.358-364).\\

It admits many properties in particular we have
\begin{itemize}
\item[i)] For all $f$ in $C_c(\mathbb{R}^d)$ we have
$$\displaystyle\int_{\mathbb{R}^d}{}^tV_k(f)(y)dy = \displaystyle\int_{\mathbb{R}^d}f(x)
\omega_k(x)dx.\eqno{(2.8)}$$ \item[ii)] For all $f$ in $
C_c(\mathbb{R}^d)$ the function $^tV_k(f)$ belongs to $
C_c(\mathbb{R}^d)$ and we have $$ \mbox{supp} f \subset
\overline{B}(o,a) \Longleftrightarrow \mbox{supp} \, {}^tV_k(f)
\subset \overline{B}(o,a), \eqno{(2.9)}$$ where
$\overline{B}(o,a)$ is the closed ball of center $o$ and radius $a
> 0$.
\item[iii)] For all $f$ in $C_c(\mathbb{R}^d)$ and $r > 0$, we
have $$\forall\; y \in \mathbb{R}^d,\; {}^tV_k(f)(ry) =
r^{2\gamma}\;{}^tV_k(f_r)(y), \mbox{ with } f_r(x) =
f(rx).\eqno{(2.10)}$$ \item[i$\nu$)] For all $a > 0$, we have
$$\forall\; y \in \mathbb{R}^d,\; {}^tV_k(e^{-a\|x\|^2})(y) =
\frac{e^{-a\|y\|^2}}{a^\gamma \pi^{d/2} c_k}.\eqno{(2.11)} $$
\end{itemize}
\hspace*{5mm} The result of the following proposition has been given
in [23] p. 363, without proof.\vspace*{3mm}\\ {\bf Proposition 2.1.}
For all $y \in \mathbb{R}^d$ we have $$supp \nu_y \cap \{x \in
\mathbb{R}^d / \|x\| = \|y\|\} \neq \emptyset. \eqno{(2.12)}$$ {\bf
Proof}\\ \hspace*{5mm}-{\bf${1}^{st}$ case:} $y \in
\mathbb{R}^{d}\backslash\{0\}$ \\Suppose to the contrary that $supp
\nu_y \cap \{x \in \mathbb{R}^d / \|x\| = \|y\|\} =\emptyset$ for
some $y $.
 Then there exists a constant $\sigma
\in ]1, + \infty[$ such that $supp \nu_y \subset \{x \in
\mathbb{R}^d / \|x\| \geq \sigma \|y\|\}$. Thus from (2.6) for all
$a > 0$ we have $${}^tV_k(e^{-a\|x\|^2})(y) = \int_{||x|| \geq
\sigma||y||}e^{-a\|x\|^2}d\nu_{y}(x).$$We put $u =
\frac{x}{\sigma} $ then  $$ {}^tV_k(e^{-a\|x\|^2})(y) =
\int_{||u|| \geq ||y||}e^{-a\sigma^{2}\|u\|^2}d\nu_{y}(u).$$ Then
$$  \; {}^tV_k(e^{-a\|x\|^2})(y) =
{}^tV_k(e^{-a\sigma^{2}\|x\|^2})(y).$$
By using the relation (2.11) we obtain
 $$
\frac{ e^{-a\|y\|^2} }{a^\gamma
\pi^{d/2}c_k}=\frac{e^{-a\sigma^{2}\|y\|^2}}{(a\sigma^{2})^\gamma
\pi^{d/2} c_k}. $$
thus
$$e^{-a\|y\|^2}=\frac{e^{-a\sigma^{2}\|y\|^2}}{\sigma^{2\gamma }}.
 $$
 If we tends  $a$ to zero we obtain $$\sigma^{2\gamma} =
1.$$ As $\gamma > 0$. Then $\sigma = 1$.
Contradiction.\\\hspace*{5mm}-{\bf${2}^{nd}$ case:} $y = 0$\\
Suppose to the contrary  that $0 \notin supp \nu_{o}$. Then there
exists $r > 0$ such that $supp \nu_{o}$ is contained in
$B^{c}(o,r)$ the complementary of the open ball $B(o,r)$ of center
$o$ and radius $r$.\\ Let $C$ be a compact contained in
$B^{c}(o,r)$ . There exists $ R > 0$ such that $$ C \subset
B(o,R). $$ We have $$\nu_{o}(C) \leq \nu_{o}(B(o,R)). $$ By
applying the relation (2.7) with $a = \frac{R}{\varepsilon}$, to
the second member of the inequality  we obtain $$\nu_{o}(C) \leq
(\frac{R}{\varepsilon})^{2\gamma}\nu_{o}(B(o,\varepsilon)).$$ Thus
$\nu_{o} = 0$. Impossible.\\ This completes the proof of the
proposition. \vspace*{3mm}\\ {\bf Theorem 2.1.} Let $(\nu_y)_{y
\in \mathbb{R}^d}$ be the family of measures defined in formula
(2.6) and let $f$ be an integrable function on $\mathbb{R}^d$ with
respect to the measure $\omega_k(x)dx$. Then for almost all $y$
(with respect to the Lebesgue measure on $\mathbb{R}^d)$, $f$ is
$\nu_y$-int\'egrable, the function $$y \mapsto \nu_y(f) =
\displaystyle\int_{\mathbb{R}^d}f(x)d\nu_y(x),$$ which will also
be denoted by ${}^tV_k(f)$, is defined almost everywhere on
$\mathbb{R}^d$ and is Lebesgue integrable. Moreover for all
bounded continuous functions $g$ on $\mathbb{R}^d$, we have the
formula
$$\displaystyle\int_{\mathbb{R}^d}{}\nu_y(f)g(y)dy =
\displaystyle\int_{\mathbb{R}^d}f(x)V_k(g)
(x)\omega_k(x)dx.\eqno{(2.13)}$$(See [7]).\newpage \noindent{\bf
Theorem 2.2.} Let $(\mu_x)_{x\in \mathbb{R}^d}$ be the family of
measures defined in formula (2.1) and let $g$ be a measurable and
bounded function on $ \mathbb{R}^d$. Then for almost all $x$ (with
respect to the Lebesgue measure on $ \mathbb{R}^d$) the function $$
x \mapsto \mu_x(g) =
\displaystyle\int_{\mathbb{R}^d}g(y)d\mu_{x}(y)$$ which also will be
denoted by $V_k(g)$, is defined almost everywhere on $
\mathbb{R}^d$, measurable and bounded. Moreover for all functions
$f$ in $C_c( \mathbb{R}^d)$ we have the formula $$
\displaystyle\int_{\mathbb{R}^d}\mu_x(g)f(x)\omega_k(x)dx =
\displaystyle\int_{\mathbb{R}^d}{}^tV_k(f)(y)g(y)dy.\eqno{(2.14)}$$
{\bf Proof}\\ \hspace*{5mm} We will divide the proof in three steps.
\\ \hspace*{5mm}i) From the properties of the operator $V_k$ we
deduce that the family of measures $(\mu_x)_{x\in \mathbb{R}^d}$
is weak-$^{*}$continuous. More precisely  for all $g$ in $C(
\mathbb{R}^d)$ the function $$ x \mapsto \mu_x(g) = V_k(g)(x) =
\displaystyle\int_{\mathbb{R}^d}g(y)d\mu_{x}(y),$$ belongs to $C(
\mathbb{R}^d)$.\\ \hspace*{5mm}ii) Let $f$ be in $C_c(
\mathbb{R}^d)$. From the relation (2.13) we deduce that for all
bounded function $g$ in $C( \mathbb{R}^d)$ we have $$
\displaystyle\int_{\mathbb{R}^d}\mu_x(g)f(x)\omega_k(x)dx =
\displaystyle\int_{\mathbb{R}^d}{}^tV_k(f)(y)g(y)dy.\eqno{(2.15)}$$\hspace*{5mm}iii)
If $g$ is a measurable and bounded function on $ \mathbb{R}^d$.
Then parts  i), ii) and Bourbaki's integration of measures Theorem
[1 , p.17] shows that the function $x \to \mu_x(g)$ exists for
almost all $x \in \mathbb{R}^d$ with respect to the Lebesgue
measure, is measurable and bounded, and the relation (2.15) is
valid for this function $g$.\vspace*{3mm}\\

 The following theorem gives the
expression of $^{t}V_k(f)$ when $f$ is radial. (See [23]).\newpage
\noindent \noindent{\bf Theorem 2.3.} For $\gamma > 0$ and for all
$f$ in $D(\mathbb{R}^d)$ radial,  we have
$$
 \forall y \in \mathbb{R}^d, \; ^{t}V_k(f)(y) = \frac{\Gamma(\gamma+ \frac{d}{2})d_k}
 {\pi^{\frac{d}{2}}\Gamma(\gamma)}\displaystyle\int_{||y||}^{+\infty}F(t)(t^2 - ||y||^2
 )^{\gamma-1}tdt,\eqno{(2.16)}$$
where $F$ is the function in $\mathcal{D}(\mathbb{R}_+) $ given by
$$f(x)= F(||x||) = F(r),\; with \,\; r = ||x||.$$
 {\bf{Examples
2.1 }}\\ \hspace*{5mm}1) When $d = 1$ and $W = \mathbb{Z}_2$, the
Dunkl intertwining operator $V_k$ is defined by ({2.1})
with for
$x \in \mathbb{R} \backslash \{0\} $ we have $d\mu_{x}(y) = {\cal
K}(x,y)dy$, where ${\cal K}$ given by the relation
(1.20).

The dual Dunkl intertwining operator $^{t}V_k$ is defined
 by (2.6)
with for all $y \in \mathbb{R}$ we have $d\nu_{y}(x) = {\cal
K}(x,y)\omega_k (x) dx$, where ${\cal K}$ and $\omega_k$  given
respectively by the relations ({1.20}) and ({1.4}).\\
 \hspace*{5mm}2) The
Dunkl intertwining operator $V_k$ of index $\gamma  =
\sum_{l=1}^{d}\alpha_l,$ $\alpha_l > 0$, associated with the
reflection group $\mathbb{Z}_2 \times ...\times \mathbb{Z}_2$ on
$\mathbb{R}^d$, is given for all $f$ in $C(\mathbb{R}^d)$ and $ x
\in \mathbb{R}^{d}_{reg} = \mathbb{R}^d \backslash \displaystyle
\bigcup_{l=1}^{d} H_l$, with $H_l = \{x \in \mathbb{R}^d/ x_l =
0\}$,  by $$ \; V_k (f)(x) = \displaystyle\int_{\mathbb{R}^d}
{\cal K}(x,y)f(y)dy \eqno{(2.17)}$$ where ${\cal K}(x,y)$ is given
by the relation (1.22). By  change of variables we obtain
$$
  \begin{array}{ccc}
    \forall \; x \in \mathbb{R}^d, \; V_k (f)(x) & = &  [
    \prod_{l=1}^d
    \frac{\Gamma(\alpha_l + \frac{1}{2})}
 {\sqrt{\pi}\Gamma(\alpha_l)}]\displaystyle\int_{[-1,1]^d}
 f(t_1x_1,t_2x_2,...,t_dx_d)\\
    &&\times\prod_{l=1}^d (1 -
    t_{l})^{\alpha_l}(1+t_l)^{\alpha_l-1}dt_1...dt_d.
  \end{array}
  \eqno{(2.18)}
  $$
(See [25] p.2964).\\ \hspace*{5mm} The dual Dunkl intertwining
operator is given for all $f$ in $C_c(\mathbb{R}^d)$ by
$$\forall\; y \in \mathbb{R}^d,\; {}^tV_k(f)(y) =
\displaystyle\int_{\mathbb{R}^d} \mathcal{K}(x,y) f(x) \omega_k
(x)dx,\eqno{(2.19)} $$ where $\mathcal{K}(x,y)$ is defined by the
relation (1.22) and $$\omega_k(x) = \prod^d_{j=1}
|x_j|^{2\alpha_j}\eqno{(2.20)}$$
\section{Absolute continuity of the representing measures of the Dunkl
 intertwining operator and of its dual}

The example 2.1 shows that when $d = 1$ and $W = \mathbb{Z}_2$ the
representing measures of the Dunkl intertwining operator and of
its dual are absolutely continuous.

 In this section we  suppose
that $d \geq 2$.

We give now the following remark which concerns the multiplicity
function $k$.\vspace*{3mm}\\ {\bf Remark 3.1.}

 Let $ \mathcal{V}' \subset \mathbb{R}^d$ be the
 $\mathbb{R}$-linear space of the subsystem
 $R' = \{\alpha \in R ; k(\alpha) \neq
 0\}$ (with $\mathcal{V}' = \{0\}$ if $R' = 0\}$
 and $\mathcal{V}'' = (\mathcal{V}')^\bot \neq
 \{0\}$. We have $\mathbb{R}^d = \mathcal{V}' \oplus
 \mathcal{V}''$. Thus all $x \in
 \mathbb{R}^d$ can be written in the form $x = x' +
 x''$ with $x' \in
 \mathcal{V}'$ and $x'' \in \mathcal{V}''$
 (see [12]). \\
 From the relations (1,9), (1,11),
 (1,16) we have
 $$\forall\; x, \lambda \in \mathbb{R}^d,\;\; K(x,\lambda)
  = e^{\langle x'',\lambda''\rangle}K(x',\lambda').$$

Using this relation and (2.13) we deduce that the measures $\mu_x$
and $\nu_y$ with $x, y \in \mathbb{R}^d$, of the integral
representations of the Dunkl intertwining operator $V_k$ and its
dual ${}^tV_k$ are of the form
$$\mu_x = \delta_{x''} \otimes \mu_{x'},\eqno{(3.1)}$$
$$\nu_y = \delta_{y''} \otimes \nu_{y'} ,\eqno{(3.2)}$$
where $\delta_{z''}$ is the Dirac measure at the point
 $z'' \in
\mathcal{V}''$.

Thus for $x, y \in \mathbb{R}^d \backslash \mathcal{V}'$, the
measures $\mu_x$ and $\nu_y$ are not absolute continuous.

 In this section we shall suppose that the multiplicity
 function $k$
 satisfies
 $$\forall\; \alpha \in R,\quad k(\alpha) > 0.\eqno{(3.3)}$$
\subsection{Absolute continuity of the measure $\mu_x$}

M.F.E. de Jeu and M. R\"osler have proved in [12] that for all $x
\in \mathbb{R}^d_{\mbox{ reg}}$ the measure $\mu_x$ of the integral
representation (2.1) of the Dunkl intertwining operator $V_k$ is
continuous

 The purpose of this subsection is to prove that for
all $x \in \mathbb{R}^d_{\mbox{reg}} $ the measure $\mu_x$ is
absolute continuous with respect to the Lebesgue measure on
$\mathbb{R}^d$, and to present some applications of this result.
\vspace*{3mm}\\ \noindent{\bf{Notations.}} We denote by \\ - $m$ the
Lebesgue measure on $ \mathbb{R}^d$.\\ - $B(\xi,r)$ the open ball of
center $\xi$ and radius $r > 0$.\vspace*{3mm}\\\hspace*{5mm}By
applying [21] p. 341 and Theorem 8.6 of [20] p.166 to the measure
$\nu_y$, $y \in \mathbb{R}^d$, we deduce that there exist a positive
function ${\cal K}^{o}(.,y)$ locally integrable on $ \mathbb{R}^d$
with respect to the Lebesgue measure, and a positive measure
$\nu_y^s$ on $ \mathbb{R}^d$ such that for every Borel set $E$ we
have $$\nu_y(E) = \displaystyle\int_{E}{\cal K}^{o}(x,y)dx +
\nu_y^s(E). \eqno{(3.4)} $$ with $${\cal K}^{o}(x,y) = \lim_{r \to
0}\frac{{\nu}_y(B(x,r))}{m(B(x,r))}.\eqno{(3.5)} $$ The measure
$\nu_y^s$ and the Lebesgue measure $m$ are mutually
singular.\vspace{5mm}\\
{\bf Remark 3.2}

When the multiplicity function satisfies
$$\forall\; \alpha \in R,\quad
 k(\alpha) = 0,\eqno{(3.6)}$$
which is equivalent to say that the subset $\mathcal{V}'$ of the
Remark 3.1 is empty, then from (3.2), (3.5)  we deduce that
$$\mathcal{K}^0(x,y) = \left\{ \begin{array}{ll}
0 &\mbox{ if } x \neq y ,\\
+ \infty &\mbox{ if } x = y.
\end{array}\right.\eqno{(3.7)}$$
Thus
$$\int_E \mathcal{K}^0(x,y) dx = 0.\eqno{(3.8)}$$
\vspace*{3mm}\\
\hspace*{5mm} For $x_{o},y \in \mathbb{R}^d$ and $n \in
\mathbb{N}^{*}$, we consider the following sequence given by
$$\overline{\Delta_n}(x_{o},y) = \sup_{0 < \frac{1}{p} <
\frac{1}{n}}\{\frac{{\nu}_y(B(x_{o},\frac{1}{p}
))}{m(B(x_{o},\frac{1}{p}))}\}.\eqno{(3.9)}$$ These functions and
their properties have been given only in the French translation [19]
of Rudin's book [20].\vspace*{3mm}\\\noindent{\bf{Lemma 3.1.}}\\
\hspace*{5mm}i) The sequence $\{\overline{\Delta_n}(x_{o},y)\}_{n
\in \mathbb{N}^{*} }$ is decreasing.\\\hspace*{5mm}ii) For $x_{o}
\in \mathbb{R}^d$ and $n \in \mathbb{N}^{*}$, the function
$\overline{\Delta_n}(x_{o},.)$ is measurable
positive.\\\hspace*{5mm}iii) For $x_{o},y \in \mathbb{R}^d$ we have
$${{\cal K}}^{o}(x_{o},y) = \lim_{n \to
\infty }\overline{\Delta_n}(x_{o},y). $$
\noindent{\bf{Proof}}\\
\hspace*{5mm} We deduce i),ii) and iii)
 from the definition of the
function $\overline{\Delta_n}(x_{o},y)$, (3.5), and the relation (3)
of [19] p.147.\vspace*{3mm}\\ \noindent{\bf{Lemma 3.2.}} For $x_0
\in \mathbb{R}^d$, $r > 0$ and for all bounded continuous function
$g$ on $\mathbb{R}^d$ we have
$$\displaystyle\int_{B(x_0,r)}V_k(g)(x)\omega_k(x)dx   =
\displaystyle\int_{\mathbb{R}^d}g(y) {\nu}_y(B(x_0,r))dy.
\eqno{(3.10)}$$\noindent{\bf{Proof}}\\ \hspace*{5mm} We deduce
(3.10) from the relation (2.13). \noindent{\bf{Proposition 3.1.}}
Let $g$ be a bounded continuous function on $ \mathbb{R}^d$. Then
for
 $x_0 \in \mathbb{R}^d$ the function ${{\cal K}}^{o}(x_0,.)$ is
integrable on $ \mathbb{R}^d$ with respect to the Lebesgue measure
and we have $$V_k(g)(x_0)\omega_k(x_0)   =
\displaystyle\int_{\mathbb{R}^d}{{\cal K}}^{o}(x_0,y)g(y)dy.
\eqno{(3.11)}$$\noindent{\bf{Proof}}\\ \hspace*{5mm} -By writing $g
=
g^+ - g^-$, we can suppose in the following that $g$ is positive.\\
\hspace*{5mm} From the relation (3.4), for $x_0 \in \mathbb{R}^d$
and $r > 0$, we have
$$\frac{1}{m(B(x_0,r))}\displaystyle\int_{B(x_0,r)}V_k(g)(x)\omega_k(x)dx
= \displaystyle\int_{\mathbb{R}^d}g(y)
\frac{{\nu}_y(B(x_0,r))}{m(B(x_0,r))}dy. \eqno{(3.12)}$$ By using
(3.9) and by applying the relation (2) of [20] p.168 to the first
member, and Fatou Lemma to the second,
 we obtain when $r$ tends to zero.
 $$\displaystyle\int_{\mathbb{R}^d}{{\cal
K}}^{o}(x_0,y)g(y)dy \leq V_k(g)(x_0)\omega_k(x_0).
 \eqno{(3.13)}$$
We replace in this inequality the function $g$ by the constant
function equal to $1$, and next we use the fact that
$$\forall \, x
\in \mathbb{R}^d, \; V_k(1)(x) = 1, $$ we deduce that
$$\displaystyle\int_{\mathbb{R}^d}{{\cal K}}^{o}(x_0,y)dy \leq
\omega_k(x_0) < +\infty.\eqno{(3.14)}$$ Then the function ${{\cal
K}}^{o}(x_0,.)$ is integrable on $ \mathbb{R}^d$ with respect to the
Lebesgue measure.\\ \hspace*{5mm} - On the other hand from the
relation (3.10) for $x_0 \in \mathbb{R}^d$ and $n,p \in
\mathbb{N}^{*}$ with $n < p$, we have
$$\frac{1}{m(B(x_0,\frac{1}{p}))}\displaystyle\int_{B(x_0,\frac{1}{p})}V_k(g)(x)\omega_k(x)dx
\leq \displaystyle\int_{\mathbb{R}^d}g(y)
\overline{\Delta}_n(x_0,y)dy.\eqno{(3.15)}$$ By using the relation
(3.13) and Lemma 3.1, we deduce that the sequence
$\{g(y)\overline{\Delta_n}(x_{o},y)\}_{n \in \mathbb{N}^{*}}$
satisfies the hypothesis of the other version of the monotonic
convergence theorem (see Theorem 5.7.27 of [21] p.234-235). By
applying this theorem to the second member of (3.15), and the
relation (2) of [20] p.168, to the first member of the same
relation, we obtain when $n$ tends to infinity
$$\displaystyle\int_{\mathbb{R}^d}{{\cal K}}^{o}(x_0,y)g(y)dy \geq
V_k(g)(x_0)\omega_k(x_0). \eqno{(3.16)}$$ We deduce (3.11) from the
relations (3.13) and
(3.16).\vspace{5mm}\\
{\bf Remark 3.3}
\begin{itemize}
\item[i)] If we suppose that the multiplicity function $k$ satisfies
(3.6), then from (1.3), (2.1) and (3.1) we have
$$\forall\; x \in \mathbb{R}^d,\;\; w_k(x) = 1,$$
and
$$V_k = Id.$$
A proof analogue to that of Proposition 3.1 and by using the
relation (3.8), we obtain for $x_0 \in \mathbb{R}^d$ and g a
positive bounded continuous function on  $\mathbb{R}^d$, the
following inequalities which are the analogue of (3.13) and (3.16) :
$$0 = \int_{\mathbb{R}^d} \mathcal{K}^0(x_0,y)g(y)dy \leq V_k(g)
(x_0)w_k(x_0),$$
$$0 = \int_{\mathcal{K}^d} \mathcal{K}^0(x_0,y) g(y)dy \geq
V_k(g)(x_0)w_k(x_0).$$ But the second member of the second
inequality is positive. Hence we obtain a contradiction.

This proof shows that the relation (3.11) is not true in this case.
\item[ii)] When the subset $\mathcal{V}'$ of the Remark 3.1 is such
that $\mathcal{V}' \neq \emptyset$, then by using (3.2) and the
result of the preceding i) we deduce that the relation (3.11) is
also not true in this case.

This remark implies that the relation (3.11) is true only under the
assumption (3.3).
\end{itemize}
{\bf{Proposition 3.2.}}\\
\hspace*{5mm} i) For all $x_0 \in \bigcup_{\alpha \in
R_+}H_\alpha$, we have for almost all $y \in \mathbb{R}^d$:
$${{\cal K}}^{o}(x_0,y) = 0.$$ \hspace*{5mm} ii) For
 $x_0 \in \mathbb{R}^d$, we have $$ \omega_k(x_0) d\mu_{x_0}(y) =
{{\cal K}}^{o}(x_0,y)dy, \; y \in \mathbb{R}^d, \eqno{(3.17)}$$
where $\mu_{x_0}$ is the measure given by the relation
(2.1).\\\noindent{\bf{Proof}}\\ \hspace*{5mm} We deduce the results
of this proposition from the relation (3.11), the fact that
$$\omega_k(x) = 0 \Longleftrightarrow x \in \bigcup_{\alpha \in
R_+}H_\alpha$$ and the properties of the measure $\mu_{x_0}$.
\vspace*{3mm}\\\noindent{\bf{Notation.}} For all $x \in
\mathbb{R}^d_{reg}$ and $y \in \mathbb{R}^d$, we put $${{\cal
K}}(x,y) = \omega_k^{-1}(x) {{\cal
K}}^{o}(x,y).\eqno{(3.18)}\vspace*{3mm}$$ \noindent{\bf{Corollary
3.1.}} The function ${{\cal K}}(.,y)$, $y \in \mathbb{R}^d$, given
by the relation (3.18) satisfies $$\forall \, x \in
\mathbb{R}^d_{reg}, \; supp {{\cal K}}(x,y) \subset
\overline{B}(0,||x||),\eqno{(3.19)} $$ where $\overline{B}(0,||x||)$
is the closed ball of center $0$ and radius
$||x||.$\vspace*{3mm}\\\noindent{\bf{Theorem 3.1.}} The representing
measure $\mu_x$ of the Dunkl intertwining operator $V_k$ satisfies
\\ \hspace*{5mm} i) For $x \in \mathbb{R}^d_{reg}$ we have
$$d\mu_x(y) = {{\cal K}}(x,y)dy, \eqno{(3.20)}$$ where ${{\cal
K}}(x,y)$ is the function given by (3.18).\\\hspace*{5mm} ii) For
all $h$ in $C(\mathbb{R}^d)$ we have $$\forall \, x \in
\mathbb{R}^d, \; \omega_k(x) V_k(h)(x) = \displaystyle\int_{
\mathbb{R}^d}{{\cal K}}^{o}(x,y)h(y)dy, \eqno{(3.21)}$$ and
$$\forall \, x \in \mathbb{R}^d_{reg}, \; V_k(h)(x) =
\displaystyle\int_{ \mathbb{R}^d}{{\cal K}}(x,y)h(y)dy.
\eqno{(3.22)}$$\noindent{\bf{Proof}}\\ \hspace*{5mm} i) The
relations (3.17) (3.18) give the result.\\
\hspace*{5mm} ii) we obtain (3.21),(3.22) from (3.17), (3.20) and
(2.1).\vspace*{3mm}\\\noindent{\bf{Remark 3.4.}}\\ \hspace*{5mm}
From this theorem and the relation (2.1) we deduce that for all $x
\in \mathbb{R}^d_{reg}$ the measure $\mu_x$ is absolutely continuous
with respect to the Lebesgue measure.\vspace*{3mm}
\\{\bf Examples 3.1}

On $\mathbb{R}^d$ with $W = \mathbb{Z}_2
 \times \cdots \times
 \mathbb{Z}_2$ the function $\mathcal{K}(x,y)$ of
 (3.20) is given
 by  (1.22) and in this case the relation
 (3.22) is true for all
 $x \in \mathbb{R}^d_{\mbox{reg}} =
 \mathbb{R}^d \backslash
  \cup^d_{\ell = 1}H_\ell$
  with $H_\ell = \{x \in \mathbb{R}^d / x_\ell = 0\}$.\vspace{5mm}\\
{\bf Corollary 3.2.}\\ \hspace*{5mm} i) For  all $x \in
\mathbb{R}^d$ and $z \in \mathbb{C}^d$ we have $$\omega_{k}(x)K(x,
-iz) = \displaystyle\int_{\mathbb{R}^d} \mathcal{K}^{o}(x,y)
e^{-i\langle y,z\rangle}dy.\eqno{(3.23)}$$\hspace*{5mm}ii) For all
$x \in \mathbb{R}^d_{\mbox{reg}}$ and $z \in \mathbb{C}^d$ we have
$$K(x, -iz) = \displaystyle\int_{\mathbb{R}^d} \mathcal{K}(x,y)
e^{-i\langle y,z\rangle}dy.\eqno{(3.24)}$${\bf Proof}

We deduce the relations (3.23),(3.24) from the relation (3) and
Theorem 3.1, ii).\vspace{5mm}

In the following proposition we give some other properties of the
function $\mathcal{K}(x,y)$.\\{\bf Proposition 3.3}
\begin{itemize}
\item[i)] For  all $x \in \mathbb{R}^d_{\mbox{reg}}$ we have
$$\displaystyle\int_{\mathbb{R}^d}\mathcal{K}(x,y)dy =
1.\eqno{(3.25)}$$
\item[ii)] For all $r > 0$, $w \in W$ and
 for  all $x \in \mathbb{R}^d_{\mbox{reg}}$ we have
 $$\mathcal{K}(wx, y) = \mathcal{K}(x,wy), a.e. y
\in \mathbb{R}^d,\eqno{(3.26)} $$ $$\mathcal{K}(rx,y) =
r^{-d}\mathcal{K}(x, \frac{y}{r}), a.e. y \in
\mathbb{R}^d.\eqno{(3.27)}$$
\end{itemize}
{\bf Proof}\\ \hspace*{5mm}We deduce these relations from Theorem
3.1 ii) and the relations\\ (3.21),(1.18).\vspace*{3mm}\\
 {\bf Corollary 3.3.}
The generalized Bessel function $J_{W}$ defined for $x \in
\mathbb{R}^{d}$ and $z \in \mathbb{C}^d$ by (see [23] p.355-356)
$$J_{W}(-ix,z) = \frac{1}{|W|}\sum_{w \in W}K(-ix,wz),
\eqno{(3.28)}$$ admits the following integral representations
$$\forall x \in \mathbb{R}^d, \; \omega_{k}(x)J_{W}(-ix,z) =
\int_{\mathbb{R}^d}E_{W}(-iz,y)
\mathcal{K}^{o}_{W}(x,y)dy,\eqno{(3.29)}
 $$ and $$\forall x \in \mathbb{R}^d_{reg},
\; \omega_{k}(x)J_{W}(-ix,z) = \int_{\mathbb{R}^d}E_{W}(-iz,y)
\mathcal{K}_{W}(x,y)dy,\eqno{(3.30)}$$ with $$E_{W}(-iz,y) =
\frac{1}{|W|}\sum_{w \in W}e^{-i\langle y,wz\rangle},
\eqno{(3.31)}$$
$$\mathcal{K}^{o}_{W}(x,y) = \frac{\omega_{k}^{-1}(x)}{|W|}\sum_{w
\in W}\mathcal{K}^{o}(wx,y),\eqno{(3.32)}$$
$$\mathcal{K}_{W}(x,y) =
\frac{1}{|W|}\sum_{w \in W}\mathcal{K} (wx,y).\eqno{(3.33)}$${\bf
Proof}\\ \hspace*{5mm} We deduce (3.29), (3.30) from the definition
of the function $J_W$ and Corollary 3.2.
 \subsection{Absolute
continuity of the measure $\nu_y$} The purpose of this subsection
is to prove that for all $y \in \mathbb{R}^d$
 the measure $\nu_y$ of the integral representation (2.6) of the dual
 Dunkl intertwining operator ${}^tV_k$ is absolute continuous with
 respect to the measure $\omega_k(x)dx$.\vspace*{3mm}\\ {\bf Lemma 3.3.}
 For $x,y_{o} \in \mathbb{R}^{d}$, we have $$
 \lim_{r \to 0}\{ \frac{{\mu}_x(B(y_0,r))}{m(B(y_0,r))}\omega_k(x) \} = {\cal
 K}^{o}(x,y_{o}),
\eqno{(3.34)}$$ where ${\cal
 K}^{o}(.,y_0)$ is the function given by the relation
 (3.4).\\
 {\bf Proof}\\
\hspace*{5mm} From Theorem 2.2 and the relation (3.11) we have
 $$\frac{{\mu}_x(B(y_0,r))}{m(B(y_0,r))}\omega_k(x) =
 \frac{1}{m(B(y_0,r))}\displaystyle\int_{B(y_0,r)}{\cal
 K}^{o}(x,y)dy.$$
 We deduce (3.31) from the relation (2) of [20] p.168.\vspace*{3mm}
 \\\hspace*{5mm}For $x,y_{o} \in \mathbb{R}^{d}$ and $n \in \mathbb{N}^{*}$, we
 consider the following sequence given by
 $$
\overline{\Lambda}_n(x,y_0) = \sup_{0 < \frac{1}{p} <
\frac{1}{n}}\{\frac{{\mu}_x(B(y_0,\frac{1}{p}))}{m(B(y_0,\frac{1}{p}))}\omega_k(x)\}.
\eqno{(3.35)}$$ {\bf Lemma 3.4.} \\\hspace*{5mm} i) The sequence
$\{\overline{\Lambda}_n(x,y_0)\}_{n \in \mathbb{N}^{*}}$ is
decreasing. \\\hspace*{5mm} ii) For $y_{o} \in \mathbb{R}^{d}$, $n
\in \mathbb{N}^{*}$ the function $\overline{\Lambda}_n(.,y_0)$ is
measurable positive. \\\hspace*{5mm} iii)For $x,y_{o} \in
\mathbb{R}^{d}$, we have $$\lim_{n \to
+\infty}\overline{\Lambda}_n(x,y_0) = {\cal K}^{o}(x,y_{o}).$$ {\bf
Proof}\\ \hspace*{5mm} We obtain  i),ii) and iii) from the
definition of $\overline{\Lambda}_n(x,y_0)$, (3.34) and the
relation (3) of [19] p.147.\vspace*{3mm}\\
 {\bf Lemma 3.5.} For $y_0 \in \mathbb{R}^d$, $r > 0$, and for all
 $f$ in $C_c( \mathbb{R}^d)$ we have $$ \displaystyle\int_{B(y_0,r)}\,^{t}V_k(f)(y) dy
 = \displaystyle\int_{
 \mathbb{R}^d}\mu_x(B(y_0,r))f(x)\omega_k(x)dx.
 \eqno{(3.36)}$$ {\bf Proof}

 We deduce (3.36) from the relation (2.14).\\\\
 {\bf Proposition 3.4.} For all $f$ in $C_c( \mathbb{R}^d)$ and
 $y_0 \in \mathbb{R}^d$ we have $$ ^{t}V_k(f)(y_0)
 = \displaystyle\int_{
 \mathbb{R}^d}{\cal K}^{o}(x,y_0)f(x)dx.\eqno{(3.37)}$$
{\bf Proof}

\hspace*{5mm} -By writing $f = f^+ - f^-$, we can suppose in the
following that $f$ is positive.\\ \hspace*{5mm} From the relation
(3.36), for $y_0 \in \mathbb{R}^d$ and $r > 0$, we have
$$\frac{1}{m(B(y_0,r))}\displaystyle\int_{B(y_0,r)}\,
^{t}V_k(f)(y)dy = \displaystyle\int_{\mathbb{R}^d}f(x)
\frac{{\mu}_x(B(y_0,r))}{m(B(y_0,r))}\omega_k(x)dx. \eqno{(3.38)}$$
By  applying the relation (2) of [20] p.168, to the first member of
(3.38), and Fatou Lemma to the second member of the same relation,
 we obtain when $r$ tends to zero. $$
 \displaystyle\int_{\mathbb{R}^d}{{\cal
K}}^{o}(x,y_0)f(x)dx \leq ^{t}V_k(f)(y_0)< +\infty. \eqno{(3.39)}$$
Thus the function ${\cal K}^{o}(.,y_0)$ is locally integrable on $
\mathbb{R}^{d}$ with respect to the Lebesgue measure.\\\hspace*{5mm}
- From (3.38)  for $y_0 \in \mathbb{R}^d$ and $n,p \in
\mathbb{N}^{*}$ with $n < p$ we have
$$\frac{1}{m(B(y_0,\frac{1}{p}))}\displaystyle\int_{
 B(y_0,\frac{1}{p})}\, ^{t}V_k(f)(y)dy \leq
  \displaystyle\int_{\mathbb{R}^d}f(x)
  \overline{\Lambda}_n(x,y_0)dx.\eqno{(3.40)}$$
  From the relation (3.39) and Lemma 3.4 we deduce
  that the
  sequence\\
   $\{f(x)\overline{\Lambda}_n(x,y_0)\} _{n \in \mathbb{N}^{*}}$ satisfies the hypothesis of the other version of the
    monotonic convergence theorem (see Theorem 5.7.27 of
     [21] p.234-235). By applying this theorem to the second member
      of (3.40),
    and the relation (2) of [20] p.168, to the first member of the same relation,
     we obtain when $n$ tends to infinity.  $$ ^{t}V_k(f)(y_0) \leq
\displaystyle\int_{\mathbb{R}^d}{{\cal K}}^{o}(x_0,y)f(x)dx.
\eqno{(3.41)}$$ We deduce (3.37) from (3.39) and (3.41).
\vspace*{3mm}\\\noindent{\bf{Theorem 3.2.}} There is a positive
function ${{\cal K}}(.,y)$, $y \in \mathbb{R}^d$, locally integrable
on $\mathbb{R}^d$ with respect to the Lebesgue measure, such that
for all $f $ in $C_c( \mathbb{R}^d)$ we have $$\forall \, y \in
\mathbb{R}^d, \; ^{t}V_k(f)(y) =
\displaystyle\int_{\mathbb{R}^d}{{\cal
K}}(x,y)f(x)\omega_k(x)dx.\eqno{(3.42)}$$\\
\noindent{\bf{Proof}}\\
\hspace*{5mm} We deduce the results from Proposition 3.4 and the
relation (3.18).\vspace*{3mm}\\ {\bf Remark 3.5}\\ \hspace*{5mm}
 Theorem 3.2 shows that for
all $y \in \mathbb{R}^d$ the measure $\nu_y$ is
absolutely
continuous with respect to the measure $\omega_k(x)dx$.
More
precisely for all $y \in \mathbb{R}^d$ we have
$$d\nu_y(x) =
\mathcal{K}(x,y) \omega_k(x)dx.$$ \noindent{\bf{Proposition 3.5.}}
For $y \in \mathbb{R}^d$ and almost $t > 0$, we have $$
\frac{1}{d_k} \displaystyle\int_{S^{d-1}} {\cal
K}(t\beta,y)\omega_k(\beta)d\sigma(\beta) =
\frac{\Gamma(\gamma+\frac{d}{2})d_k}{\pi^{\frac{d}{2}}\Gamma(\gamma)}
t^{2-2\gamma-d}(t^2 - ||y||^2)^{\gamma-1}1_{]||y||,+\infty[}(t),
\eqno{(3.43)}$$ where $1_{]||y||,+\infty[}$ is the characteristic
function of the interval ${]||y||,+\infty[}$.\vspace*{3mm}\\
\noindent{\bf{Proof.}}\\ \hspace*{5mm}Let $f$ be a radial  function
in ${\cal D}( \mathbb{R}^d).$ From (3.42) we have  $$\forall \, y \,
\in \mathbb{R}^d, \; \,^{t}V_k(f)(y) =
\displaystyle\int_{\mathbb{R}^d} {\cal
K}(x,y)F(||x||)\omega_k(x)dx,$$ where $F$ is the function in ${\cal
D}( \mathbb{R_+})$ given by $$\forall \, x \in \mathbb{R}^d, \; f(x)
= F(||x||).$$ Using (1.6) and the fact that the support of ${\cal
K}(x,y)$ is contained in the set $\{x \in \mathbb{R}^d/ ||x|| \geq
||y||\}$, we obtain $$\forall \, y \, \in \mathbb{R}^d, \;
\,^{t}V_k(f)(y) = \displaystyle\int_{||y||}^{+\infty}(
\displaystyle\int_{{S}^{d-1}} {\cal
K}(t\beta,y)\omega_k(\beta)d\sigma(\beta))F(t)t^{2\gamma+d-1}dt.$$
By applying Theorem 2.3  we deduce that for almost all $t >
0$:$$t^{2\gamma+d-1}\displaystyle\int_{{S}^{d-1}} {\cal
K}(t\beta,y)\omega_k(\beta)d\sigma(\beta)) =
\frac{\Gamma(\gamma+\frac{d}{2})d_k}{\pi^{\frac{d}{2}}\Gamma(\gamma)}
t^{2-2\gamma-d}(t^2 - ||y||^2)^{\gamma-1}1_{]||y||,+\infty[}(t).$$
We obtain (3.43) from this relation.\vspace*{3mm}\\
{\bf Example
3.2}

From the relations (2.17), (2.18) we deduce that the relation (3.42)
of the dual Dunkl intertwining operator ${}^tV_k$ on $\mathbb{R}^d$
with the reflection group $\mathbb{Z}_2 \times \mathbb{Z}_2 \times
... \times \mathbb{Z}_2$ can also be written for all
 $f$ in
$C_c(\mathbb{R}^d)$ in the form
\begin{eqnarray*}
\forall\; y \in \mathbb{R}^d, {}^tV_k (f)(y) &=&
\Big[\prod^d_{j=1} \frac{\Gamma(\alpha_j +
\frac{1}{2})}{\sqrt{\pi}\Gamma(\alpha_j)}\Big]
\displaystyle\int_{|x_1|> |y_1|}... \displaystyle\int_{|x_d|>
|y_d|}f(x_1,...,x_d) \\ &&\times \Big[ \prod^d_{j=1} (|x_j| -
y_j)^{\alpha_j}(|x_j| + y_j)^{\alpha_j -1} \Big] dx_1 ...dx_d .
\end{eqnarray*}
\section{Applications}

In this section we suppose that the multiplicity function $k$
satisfies the assumption (3.3).
\subsection{First application}
{\bf Theorem 4.1.} We have  $$\forall \, x \in \mathbb{R}^d,\;
\lim_{\|z\|\rightarrow + \infty}\{\omega_{k}(x) K(-ix,z)\} = 0.
\eqno{(4.1)}$$ and $$\forall \, x \in \mathbb{R}^d_{reg},\;
\lim_{\|z\|\rightarrow + \infty} K(-ix,z) = 0. \eqno{(4.2)}$$
{\bf
Proof}

- From Corollary 3.2 i), for  all $x \in \mathbb{R}^d$, and
 $z \in \mathbb{R}^d$ we have $$\omega_{k}(x)K(-ix,z) =
\displaystyle\int_{\mathbb{R}^d}\mathcal{K}^{o}(x,y) e^{-i\langle
y,z\rangle}dy.$$ As for  $x \in \mathbb{R}^d$ the function
$\mathcal{K}^{o}(x,.)$ is integrable with respect to the Lebesgue
measure on $\mathbb{R}^d$, then we obtain the relation (4.1) from
Riemann-Lebesgue Lemma for the usual Fourier transform on
$\mathbb{R}^d$.\\ \hspace*{5mm}- Corollary 3.2 ii) and the same
proof give the relation (4.2).\vspace*{3mm}\\ {\bf Remark 4.1.}

Let $C$ denotes the Weyl chamber attached with the  positive
subsystem $R_+$, $$C = \{x \in \mathbb{R}^d / \langle \alpha,
x\rangle > 0, \mbox{ for all } \alpha \in R_+\},$$ and for $\delta
> 0$, $$C_\delta = \{x \in C/ \langle \alpha,x\rangle > \delta\|x\|, \mbox{ for all }
\alpha \in R_+\}.$$ M.F.E de Jeu and M. R\"osler   have proved in
[12] the following behaviour for the Dunkl kernel $K(x, -iz)$,
uniform for the variable tending to infinity in cones $C_\delta$ :
There exists a constant non-zero vector $v = \{v_w\}_{w \in W}$ such
that for all $x \in C, w \in W$ and each $\delta > 0$,
$$\lim_{\|z\| \rightarrow + \infty, z \in C_\delta}
\sqrt{\omega_k(x) \omega_k(z)}e^{i(wx,z)}K(wx, -iz) = v_w.
\eqno{(4.3)}$$ {\bf Corollary 4.1.} We have  $$\forall \, x \in
\mathbb{R}^d,\; \lim_{\|z\|\rightarrow + \infty}\{\omega_{k}(x)
J_{W}(-ix,z)\} = 0. \eqno{(4.4)}$$ and $$\forall \, x \in
\mathbb{R}^d_{reg},\; \lim_{\|z\|\rightarrow + \infty}
J_{W}(-ix,z) = 0. \eqno{(4.5)}$$ {\bf Proof} \\ \hspace*{5mm} We
deduce (4.4),(4.5) from Corollary 3.3 and Theorem 4.1.
\subsection{Second application}
 {\bf Theorem 4.2.}
 For  all function $h$ in $C(\mathbb{R}^d)$ we have
 $$\begin{array}{ll}
\forall\; t > 0,\;
\displaystyle\displaystyle\int_{S^{d-1}}V_k(h)(t\xi)\omega_k(\xi)d\sigma(\xi)
= \frac{\Gamma(\gamma + \frac{d}{2})d_k}{\pi^{d/2}\Gamma(\gamma)}
t^{2 - 2\gamma - d}\times\\ \hspace*{8cm}
\displaystyle\displaystyle\displaystyle\int_{B(0,t)} h(y)(t^2 -
\|y\|^2)^{\gamma-1}dy,
\end{array}\eqno{(4.6)}$$
where $B(0,t)$ is the open ball of center $0$ and radius
$t$.\vspace*{3mm}\\ {\bf Proof}

\hspace*{5mm} From Theorem 3.1 ii) the relation (3.43) and Fubini's
theorem, for almost all $t \in ]0,+\infty[$ we have
$$\begin{array}{lll}\displaystyle\displaystyle\int_{S^{d-1}}V_k
(h)(t\xi)\omega_k(\xi)d\sigma(\xi) &=&
\displaystyle\int_{R^{d}}[\displaystyle\int_{S^{d-1}}\mathcal{K}(t\xi,y)d\sigma(\xi)]h(y)dy
\\\\&=&ct^{2-2\gamma-d}\displaystyle\int_{R^{d}}h(y)(t^2 -
||y||^2)^{\gamma-1}1_{]||y||,+\infty[}(t)dy\end{array}$$ where $$c
=
\frac{\Gamma(\gamma+\frac{d}{2})d_k}{\pi^{\frac{d}{2}}\Gamma(\gamma)}.$$
By using the spherical coordinates we obtain
$$\begin{array}{lll}\displaystyle\displaystyle\int_{S^{d-1}}V_k
(h)(t\xi)\omega_k(\xi)d\sigma(\xi)
=ct^{2-2\gamma-d}\displaystyle\int_{0}^t(\displaystyle\int_{S^{d-1}}h(\varrho
\eta)d\sigma(\eta))(t^2 -
\varrho^2)^{\gamma-1}\varrho^{d-1}d\varrho.\end{array}$$ Thus
$$\begin{array}{lll}
\displaystyle\displaystyle\int_{S^{d-1}}V_k(h)(t\xi)\omega_k(\xi)d\sigma(\xi)
= c t^{2-2\gamma-d}
\displaystyle\displaystyle\displaystyle\int_{B(0,t)} h(y)(t^2 -
\|y\|^2)^{\gamma-1}dy.
\end{array}$$
Using the properties of the operator $V_k$ we deduce that the
first member of this relation is continuous on $]0,+\infty[$. The
second member possesses also the same property. Then this relation
is true for all $t \in ]0,+\infty[$.\vspace*{3mm}\\ {\bf Remarks
4.2.}

i) \hspace*{5mm} From the relations (2.1),(1.20) we deduce that
for  $d= 1$ the analogue of the relation (4.6) is of the form
$$\forall \, x > 0, \; \frac{V_k(h)(x) + V_k(h)(-x) }{2} =
\frac{\Gamma(\gamma+\frac{1}{2})}{\pi^{\frac{1}{2}}\Gamma(\gamma)}x^{1-2\gamma}
\displaystyle\int_{-x}^{x}h(y)(x^2-y^2)^{\gamma-1}dy.$$

ii) By using the $\omega_k$-harmonic polynomials
Y.Xu has proved in
[25] the relation (4.6) for $t = 1$ and $h$ a $\omega_k$-harmonic
polynomial.

\subsection{Third application}\hspace*{5mm}
\hspace*{5mm}  The  generalized (or Dunkl) translation operators
$\tau_x$, $x \in \mathbb{R}^d$, are  defined    on ${\cal
E}(\mathbb{R}^d)$ by $$
\forall \, y \in \mathbb{R}^d, \;\tau_{x}f(y)= (V_k)_x (V_k)_y[(V_k)^{-1}(f)(x+y)]. %
$$ (See [24] p.33-35).\\ \hspace*{5mm}They satisfies many
properties, in particular we have $$\tau_{x}f(o) =f(x), \;
\tau_{x}f(y) = \tau_{y}f(x), \;\tau_{x}(1)(y) = 1.
\eqno{(4.7)}$$\hspace*{5mm} At the moment an explicit formula for
the generalized (or Dunkl)
 translation operators is known only in the following
 cases( see [18] [22]). \\
{\bf{1$^{st}$ cas }}: $d = 1$ and $W = \mathbb{Z}_2$. \\
\hspace*{5mm}For all $f$ in $C(\mathbb{R})$ and $ y \in
\mathbb{R}$ we have $$\begin{array}{ccc}
  \tau_{x}f(y) & = &  \frac{1}{2}\displaystyle\int_{-1}^{1}f(\sqrt{x^2 + y^2 -2xyt})
  (1+\frac{x-y}{\sqrt{x^2 + y^2 -2xyt}})\Phi_k(t)dt\\
   & + &  \frac{1}{2}\displaystyle\int_{-1}^{1}f(-\sqrt{x^2 + y^2 -2xyt})
  (1-\frac{x-y}{\sqrt{x^2 + y^2 -2xyt}})\Phi_k(t)dt,
\end{array}\eqno{(4.8)}$$ where $$\Phi_k(t) =
\frac{\Gamma(k+\frac{1}{2})}{\sqrt{\pi}%
\Gamma(k)} (1+t)(1-t^2)^{k-1}.$$ If we consider $f$ even and we
make the change of variables  $u = yt$, we obtain for all $ y \in
\mathbb{R}\backslash\{o\}$: $$ \tau_{x}f(y) =
\frac{\Gamma(k+\frac{1}{2})}{\sqrt{\pi}\Gamma(k)}|y|^{-2k}
\displaystyle\int_{-|y|}^{|y|}F(\sqrt{x^2 + y^2 -2xu})
  (|y|-u)^k(|y|+u)^{k-1}du,\eqno{(4.9)}$$ where $F$ is the
  restriction of $f$ on $[0,+\infty[$.\\ By using the relations
  (3.22),(1.22) we deduce that $$\forall \, y
\in \mathbb{R}\backslash\{o\}, \tau_{x}f(y)  =
 V_k[F(\sqrt{x^2 +
y^2 -2x.})](y).\eqno{(4.10)}$$
 {\bf{2$^{nd}$ cas:}} $d \geq 2$\\ \hspace*{5mm} For all
$f$ in ${\cal E}(\mathbb{R}^d)$ radial, we have $$\forall \, y \in
\mathbb{R}^d, \;  \tau_{x}f(y)  =
\int_{\mathbb{R}^d}F(\sqrt{||x||^2 + ||y||^2 -2\langle
x,\eta\rangle})d\mu_{x}(\eta),\eqno{(4.11)}$$ where $F$ is the
function on $[0,+\infty[$ given by $f(x) = F(||x||)$.\\
\hspace*{5mm} For $f$ in ${\cal E}(\mathbb{R}^d)$, even for $d = 1$
and radial for $d \geq 2$, the relations (4.10),(4.11) and (3.20)
implies that for all $x \in \mathbb{R}^d$ we have $$\forall \, y \in
\mathbb{R}^d_{reg}, \;  \tau_{x}f(y)  =
\int_{\mathbb{R}^d}F(\sqrt{||x||^2 + ||y||^2 -2\langle
x,\eta\rangle}){\cal K}(y,\eta)d\eta.\eqno{(4.12)}$$
 {\bf Theorem 4.3.} For all $x \in
 \mathbb{R}^d\backslash\{o\}$
 and $y \in \mathbb{R}^d_{reg}$, there
  exists a positive function
 ${\cal W}(x,y,(t,t'))$ satisfying
 $$\int_{[0,+\infty[\times\mathbb{R}^{d-1}}{\cal W}
 (x,y,(t,t'))dtdt' = 1,$$
 such that for all $f$ in ${\cal E}(\mathbb{R}^d)$,
 even for $d = 1$
and radial for $d \geq 2$,  we have $$\tau_{x}f(y)  =
\int_{[0,+\infty[\times \mathbb{R}^{d-1}}F(t){\cal
W}(x,y,(t,t'))dtd{t'}.\eqno{(4.13)}$$
\noindent{\bf{Proof}}\\
\hspace*{5mm} We deduce the results of this theorem from the
relation (4.12), the change of variables: $$t = \sqrt{||x||^2 +
||y||^2 -2\langle x,\eta\rangle},\, t'_1 = \eta_2, ..., t'_{d-1} =
\eta_d,$$ and the properties (4.7).\vspace*{3mm}\\
 {\bf Remark 4.3.}\\
\hspace*{5mm}Theorem 4.3 shows that the measure $\varrho_{x,y}^k$(
see Theorem 5.1 of [18]) of the integral representation of the
generalized (or Dunkl) translation operators $\tau_x$, $x \in
\mathbb{R}^d\backslash\{o\}$ corresponding to the preceding cases,
are absolute continuous with respect to the Lebesgue measure.


\begin{thebibliography}{4}
\bibitem{B}{\bf N.BOURBAKI.}{\em \'El\'ements de Math\'ematique.
Fascicule XXI. Livre V. Int\'egration-Chapitre V-Int\'egration des
mesures. Hermann-Paris 1967.}

\bibitem{Ch}
{\bf L.CHEREDNIK.}{\em A unification of the Knizhnik-Zamolodchikov
equations and Dunkl operators via affine Hecke algebras. Invent.
Math. 106, 1991, p. 411-432.}

\bibitem{DI}{\bf J.F.van DIEJEN.}{\em Confluent hypergeometric orthogonal polynomials
related to the rational quantum Calogero system with harmonic
confinement. Comm.Math.Phys.188, 1997, p 467-497.}

\bibitem{D1}{\bf C. F. DUNKL.  }
{\em Differential-difference operators associated to reflection
groups. Trans. Am. Math.
 Soc. 311, 1989, p. 167 - 183. }
\bibitem{D2}{\bf C.F. DUNKL.  } {\em Integral kernels with reflection group invariance.
Can . J . Math. 43, 1991, p. 1213-1227. }

\bibitem{D3}{\bf C.F.DUNKL .  }{\em Hankel
transforms associated  to finite
reflection groups
 Contemp. Math.  138, 1992, p. 123-138 .}

\bibitem{HA}{\bf L. GALLARDO and K. TRIM\`ECHE}
 {\em Un analogue d'un th\'eor\`eme de Hardy pour la transformation
 de Dunkl.C.R. Acad. Sci. Paris, S\'erie I, t. 334, 2002, p. 849 -
 854}.
\bibitem{Hu}{\bf G.J.HECKMAN.   }{\em An elementary approach to the
 hypergeometric
shift operators of Opdam. Invent.Math. 103, 1991, p. 341-350.}

\bibitem{H}{\bf J.E.HUMPHREYS.   }{\em Reflection groups and Coexter groups. Cambridge
Univ. Press. Cambridge, England, 1990.  }

\bibitem{HI}{\bf K.HIKAMI. }{\em Dunkl operators formalism for quantum many-body problems associated
with classical root systems.J.Phys.Soc.Japan, 65, 1996, p.
394-401.}

\bibitem{J}{\bf M.F.E.de JEU.  }{\em The Dunkl transform.
Invent.Math. 113 , 1993, p. 147-162.}

\bibitem{r3}{\bf M.F.E de Jeu and M. R\"OSLER.}
{\em Asymptotic analysis for the Dunkl kernel.J.Aprox. Theory 119,
2002, p.110-126. }

\bibitem{K}{\bf S.KAKEI}.{\em Common algebraic structure for the
 Calogero-Sutherland models.
J.Phys. A 29, 1996, p. 619-624.}

\bibitem{LV}{\bf L. LAPOINTE and L.VINET.}{\em Exact operator solution of
the
 Calogero-Sutherland model.
Comm.Math.Phys.178, 1996, p.425-452.}

\bibitem{M}{\bf M.L.MEHTA.   }{\em
Random matrices and statistical theory of energy levels. Academic
Press, New York, 1967.}
\bibitem{O}{\bf E.M.OPDAM.   }{\em Harmonic analysis for certain
 representations of graded Hecke algebras. Acta Math.175, 1995, p.
 75-121.}

\bibitem{Ros}{\bf M. R\"OSLER.}
{\em Positivity of  Dunkl's intertwining Operator. Duke Math. J.
98, 1999, p.445-463. }

\bibitem{Ros}{\bf M. R\"OSLER.}
{\em A positive radial product formula for the Dunkl kernel.
Trans. Amer. Math. Soc., 355, 2003, p.2413-2438.}

\bibitem{Rud1}{\bf W. RUDIN.}
{\em Analyse r\'eelle et complexe. Masson et C$^{ie}$, \'Editeurs.
1975. }

\bibitem{Rud2}{\bf W. RUDIN.}
{\em Real and complex analysis. Mc Graw Hill Inc. Second edition.
1966,1974. }

\bibitem{S}{\bf L. SCHWARTZ.}
{\em Analyse III, calcul int\'egrale. Hermann- \'Editeurs des
Sciences et des Arts. 1993. }

\bibitem{karni4}{\bf S.\ THANGAVELU  and Y.XU.\ }%
Convolution operator and maximal functions for Dunkl transform.
Arxiv math CA/ 0403049 2004.

\bibitem{T4}
{\bf K.TRIM\`ECHE.}{\em The Dunkl intertwining operator on spaces
of functions and distributions and integral representation of its
dual.  Integ. Transf. and Special Funct. Vol. 12, $N^o$4, 2001, p.
349-374 }.

\bibitem{T5}
{\bf K.TRIM\`ECHE. }{\em Paley-Wiener theorems for Dunkl transform
and Dunkl translation operators.Integ. Transf. and Special Funct
Vol.13, 2002, p.17-38.}

\bibitem{X}{\bf Y.XU.} {\em Integration of intertwining operator for $h$-harmonic
polynomials associated to reflection groups. Proc. Amer. Math.
Soc. Vol. 125, $N^o$ 10, 1997, p. 2963 - 2973}.
\end{thebibliography}
\end{document}